\documentstyle[11pt]{article}
\begin{document}
\makeatletter \@addtoreset{equation}{section} \makeatother
\renewcommand{\theequation}{\thesection.\arabic{equation}}
\baselineskip 15pt

\title{\bf Book Review: \\ ``A Beautiful Math. John Nash, Game Theory, and 
the Modern Quest for a Code of Nature'' \\
by Tom Siegfried \\ {\small\bf Joseph Henry Press
(2006)}}
\author{Paul B. Slater\footnote{e-mail: slater@kitp.ucsb.edu}\\
{\small ISBER,} \\
{\small University of California, Santa Barbara, 93106.}}
\date{}
\maketitle

In the collection of modeling paradigms available to scientific researchers, none has proved perhaps as productive and insightful as the {\it game-theoretic} one--in the development of which, the eminent mathematician John Forbes Nash played a most pivotal role (having introduced the concept of {\it non-cooperative} games [MR0043432, MR0031701, MR0035977], complementing the landmark  [zero-sum, minimax] work of Von Neumann and Morgenstern [MR0011937 and MR2316805]).  This is the interestingly and certainly well-defended and documented--although not universally accepted--position taken by Tom Siegfried in his concise book, "The Beautiful Mind: John Nash, Game Theory, and the Modern Quest for a Code of Nature".  (The distinguished economist and {\it game-theorist} Ariel Rubinstein is a notable nay-sayer--"the challenges facing the world today are far too complex to be captured by any matrix game" (p. 68). Further, the Nash equilibrium is a static concept, and the dynamics required to achieve it are problematical, with pay-offs being assumed to depend  {\it continuously} on the various decisions taken  [I. Ekeland, "Game theory: Agreeing on strategies", Nature 400, 623-624 (12 August 1999)]. While Siegfried's book is essentially concerned with the application of equilibrium concepts, the author does note (p. 238) a reviewer's comments that, in some cases, a more appropriate description might be provided by a chaotic system [cf. MR2169537, MR2033167, MR1963719].)

Tom Siegfried is a well-respected, often-honored science journalist (science editor of the Dallas Morning News for many years) and the author of two previous popular science books  ("The bit and the pendulum: from quantum computing to M theory-- the new physics of information"  [Wiley 2000] and  "Strange matters: undiscovered ideas at the frontiers of space and time"  [Joseph Henry Press 2002]).  

The main title of Siegfried's latest effort is an obvious play on that of the critically successful "A Beautiful Mind" of Silvia Nasar (MR1631630, MR1646747, MR1927043), which examined  in remarkable depth and detail the brilliant, but troubled career of Nash. (He did not cooperate in the preparation of his biography.) The cinematic version garnered four major Academy Awards and four Golden Globe Awards, but was viewed somewhat skeptically, it would seem, by much of the scientific community--though the widespread attention drawn 
to the "romance" of mathematics and mathematicians was, by and large, not unwelcome (MR2927043). The algebraic geometer David Bayer served as a mathematical consultant on the movie. 

It might be noted that--despite the main title--there is little mathematics {\it per se} discussed in Siegfried's book, and certainly the author's stress is on the widespread usefulness and applicability of the mathematical results of Nash, rather than any elegance or "beauty" they may exhibit. (No equations are given, perhaps in accordance with the advice given to Stephen Hawking by his publisher that $n$ equations reduces the readership of a popular science treatment by $2^{-n}$. If a narrower, simply academic audience had been targeted, "Uses of Game Theory in Sociophysics" might have been an apt title to describe the subject matter covered.)

After the appearance of the book, the 2007 Nobel Memorial Prize in Economic Sciences was
awarded equally to Leonid Hurwicz, Eric Maskin and Roger Myerson "for having laid the foundations of mechanism design theory", a branch of game theory that allows people to distinguish those situations in which markets work well from those in which they do not.
(Let us remark that one possible important topical application of [quantum] game-theoretic analysis might be to the design of international treaties pertaining to climate-change issues [MR2383839, MR2245475, MR2239192 and MR2180222]. Myerson, in particular, has emphasized Nash's work "as one of the outstanding intellectual advances of the twentieth century" [pp. 51-52].) John Nash shared the 1994 Prize with Richard Selten and John Harsanyi  "for their pioneering analysis of equilibria in the theory of non-cooperative games", while the 2005 Prize--as discussed by Siegfried [pp. 70-71]--was given to Robert Aumann and Thomas Schelling "for having enhanced our understanding of conflict and cooperation through game theory".

The eminent topologist and geometer John Milnor has suggested that the subsequent mathematical work of Nash, pertaining to embedding theorems and partial differential equations, among other topics, was "far more rich and important"--at least from a strictly mathematical viewpoint--than the "ingenious but not surprising application of well-known (fixed-point) methods" (von Neumann famously called the application "trivial") used for the game-theoretic results. But to the extent a "new work increases[s] our understanding of the real world,...Nash's thesis was nothing short of revolutionary" (MR1646747). Siegfried--for simplicity and relevance sake, evidently--only mentions the existence of such later, much-esteemed work of Nash, and just briefly discusses (pp. 58-59) the earlier game-theoretic-related use of fixed-point methods.

Siegfried ambitiously (and perhaps with some admirable, impressive hubris) seeks to complement the biographical work of Nasar, only quickly touching upon Nash's personal 
life, while trying to place (still in a popular fashion) his work on non-cooperative games in the broad development (from Roman times onwards) of human thought across many disciplines.  
Following the Preface and Introduction, the numbered chapters are: (1) (Adam) "Smith's Hand: Searching for the Code of Nature"; (2) "Von Neumann's Games: Game theory's origins"; (3) "Nash's Equilibrium: Game theory's foundations"; (4) (John Maynard) "Smith's Strategies: Evolution, altruism and cooperation"; (5) "Freud's Dream: Games and the brain"--concerned with the emerging important fields of neuroeconomics and {\it behavioral} game theory; (6) "Seldon's Solution: Game theory, culture and human nature"; (7) "Quetelet's Statistics and Maxwell's Molecules: Statistics and society, statistics and physics"; (8) (the actor Kevin, {\it not} Sir Francis) "Bacon's Links: Networks, society, and games" (for a critique of some of the work discussed here, in particular, "scale-free networks" and power laws, see MR2241756); (9) "Asimov's Vision: Psychohistory, or sociophysics"; (10) "Meyer's Penny: Quantum fun and games"; and (11) "Pascal's Wager: Games, probability, information, and ignorance".  There is also an "Epilogue", an "Appendix: Calculating a Nash Equilibrium" and a "Further Reading" guide, plus a highly-extensive scholarly set of end notes. Though much in the nature of independent essays, chapters conclude with paragraphs indicating the transition to new--but connected, it is contended--subject matter in the following chapters. (Certainly, there are far too many specific Nash-equilibrium applications for the author to have mentioned, but let me simply indicate an additional one, that to the traffic network equilibrium problem [Patriksson, M. (1994) The Traffic Assignment Problem - Models and Methods, VSP, Utrecht].)

Frequent reference --as indicated in the book's subtitle--is given to the old Roman notion of a "Code of Nature" ({\it Jus Naturale})--a theory that posits the existence of a law whose content is set by nature and that therefore has validity globally. Also, somewhat entertainingly and with a light touch, Siegfried occasionally refers to non-strictly scholarly literature--that is, popular fiction and, in particular, Isaac Asimov's epic science-fiction "Foundation Trilogy". (The title of Chapter 6 pertains to the hero of the Trilogy, a mathematician Hari Seldon, who created a community of scientists dedicated to manipulating the future. Siegfried returns many times to him and his objectives in creating a system, "psychohistory", to predict trends, governmental downfalls, wars,...)

The author's enthusiastic and intelligent, well-researched (numerous interviews and preprint references) treatment of many cutting-edge scientific topics is certainly a significant, outstanding  contribution to the popular science literature. Whether or not the diverse subjects covered are highly appropriately unified under the game-theoretic banner, and directly attributable to the influence of John Nash, in line with the main thesis of this book, surely remains an open matter for considerable serious and lively discussion among mathematicians, physicists, economists, social scientists, biologists,...Nash  is, beyond doubt, a brilliant, highly original mathematician. His career and oeuvre has not been that of a profound social philosopher, on the order of Adam Smith, as a casual peruser of the book--with its interesting, repeated mention of the "Modern Quest for a Code of Nature" and of Nash himself and Nash equilibrium--may, at first, be led to think. He clearly laid the foundation for many additional, important developments, but did not directly participate in them (if for no other reasons than those of the well-publicized  personal nature).

\end{document}